\documentclass{elsarticle}

\usepackage{lineno,hyperref}
\usepackage{a4wide}
\usepackage{amsmath}
\usepackage{amsthm}
\usepackage{amssymb}
\usepackage{amsfonts}
\usepackage{mathrsfs}
\usepackage{stmaryrd}
\usepackage{booktabs}
\usepackage{wrapfig}
\usepackage{graphicx}
\usepackage{epstopdf}
\usepackage{paralist}
\usepackage{subcaption}
\usepackage{xcolor}
\usepackage{enumitem}
\usepackage{siunitx}
\usepackage{multirow}
\usepackage{xspace}
\usepackage[capitalise]{cleveref}
\usepackage{todonotes}
\usepackage{bm}
\usepackage{nicefrac}
\usepackage[labelfont=bf,margin=1in,labelsep=endash,justification=raggedright]{caption}

\newcommand{\eg}{e.g. }

\theoremstyle{plain}

\newcommand{\Dumux}{{Du\-Mu$^\text{x}$ }}

\makeatletter
\def\ps@pprintTitle{%
	\let\@oddhead\@empty
	\let\@evenhead\@empty
	\def\@oddfoot{}%
	\let\@evenfoot\@oddfoot}
\makeatother

\begin{document}
	
	\begin{frontmatter}
		
		\title{A \Dumux Framework for Data-Driven Multi-Scale Parametrizations}
		
		\author[addressIWS]{Edward Coltman \corref{mycorrespondingauthor}}
		\ead{edward.coltman@iws.uni-stuttgart.de}
		\author[addressIWS]{Martin Schneider}
		\ead{martin.schneider@iws.uni-stuttgart.de}
		\author[addressIWS]{Rainer Helmig}
		\ead{rainer.helmig@iws.uni-stuttgart.de}
		
		\cortext[mycorrespondingauthor]{Corresponding author}
		
		\address[addressIWS]{Institute for Modelling Hydraulic and Environmental Systems,
			University of Stuttgart,
			Pfaffenwaldring 61,
			70569 Stuttgart, Germany}
		
		\begin{abstract}
			
			Presented in this work is a framework for the data-driven determination of multi-scale porous media parametrizations. Simulations of flow and transport in a porous medium at the REV scale, although efficient, require well defined parameters that represent pore-scale phenomena to maintain their accuracy. Determining the optimal parameters for this often require expensive pore-scale calculations. This work outlines a series of four steps where these parameters can be calculated from pore scale data, their solutions generalized with a convolutional neural network, and their content better understood with descriptive pore metrics.
			 
		\end{abstract}
		
	\end{frontmatter}
	
	%%%%%%%%%%%%%%%%%%%%%%%%%%%%%%%%%%%
	%%%%%%%%%%%%%%%%%%%%%%%%%%%%%%%%%%%
	%%%%%%%%%%%%%%%%%%%%%%%%%%%%%%%%%%%
	\section{Introduction}
	\label{sec:introduction}
	
	Evaluations on the Representative Elementary Volume (REV) scale are ubiquitous in the field of environmental modelling and porous medium flow research \cite{bear1988a}. Although these models are efficient and often very accurate, phenomena seen on higher resolution scales must be ignored on this averaged scale. In some more complex systems, these pore-scale phenomena can play a significant role in the dynamics of the system. While neglecting these effects can worsen model results, a full evaluation on the pore-scale of these complex systems is often not realistically possible. For this reason, additional multi-scale closure terms can be added to the averaged scale models to represent the pore-scale dynamics on the REV scale. The use of machine learning tools is becoming more common within this field \cite{wang2021a,tahmasebi2020a}, specifically for efficiently replacing any multi-scale closure problems with data-driven models trained on pore-scale solutions \cite{taghizadeh2022a}. Considerable work has been done investigating the permeability parameter using this method\cite{gaerttner2021a,mostaghimi2006a}, and further work has developed evaluating transport and multi-phase flow parameters \cite{lasseux2021a}.
	
	In this work, and the corresponding numerics environment, a framework is outlined for the data-driven parametrizations of these terms, as well as methods for evaluating these parametrizations with pore-scale metrics. The framework is divided into 4 parts: Pore-Scale Data Development (Section \ref{sec:porescaledata}), Averaging and REV scale model parametrization (Section \ref{sec:averagescale}), parameterization generalization with machine learning tools (Section \ref{sec:mlgeneralization}), and analysis using descriptive metrics (Section \ref{sec:metrics}). A brief summary and outlook are provided in Section \ref{sec:conclusions}. All code and code examples are provided in \cite{PubModule}.

	%%%%%%%%%%%%%%%%%%%%%%%%%%%%%%%%%%%
	%%%%%%%%%%%%%%%%%%%%%%%%%%%%%%%%%%%
	
	\section{Pore-Scale Evaluation and Data Production} \label{sec:porescaledata}
	
	Paramount to the success of any data-driven modelling effort is the production and curation of quality data. In this first step, training data, in the form of high resolution pore-scale simulations, is to be created and logged. 
	
	\subsection{Pore Scale Model}
	Models at this detailed pore-scale should be chosen such that the target phenomena be captured. For example, when evaluating flow parameters such as permeability of a Forchheimer coefficient, direct numerical simulations of the flow field should be performed with models such as the Navier Stokes equations or a Lattice-Boltzman model. If transport properties are to be evaluated, transport models should be included, and if multi-phase parameters are to be evaluated, multiphase extensions, such as phase field models or volume of fluids methods, to the direct numerical simulation must be made. 
	
	In this framework, the direct numerical simulation of pore-scale flow and transport is developed using \Dumux 's freeflow simulation environment \cite{koch2020a}. 
		
	\subsection{Pore Geometries}	
	
	Another key challenge of pore-scale analysis is the definition of a pore geometry upon which the simulation domain can be laid. With recent developments in micro-CT scanning technology, the determination of real pore geometries has become more feasible \cite{prodanovic2015a}. Within this framework, the import and export of pore geometries to a numerical grid is done using DUNE-SubGrid \cite{graeser2009a}, with a grid manager developed to read binary or raster images. Although this interface can use images of any origin, additional generators are provided for the development of simple as well as random pseudo pore-geometries. Each generator provided defaults to a standard 2D resolution of 200x200 pixels. 
		
	\subsubsection{Pre-processing}
	In order to ensure their usability, each unit cell is checked for connectivity using a graphical component traverse, further described in Section \ref{subsubsec:connectivity}. Any disconnected pore space is then excluded directly by altering the image. In the case of high discontinuity, the pore space should be redeveloped before simulation. In addition, the periodic connectivity of these shapes is evaluated to be used on periodic simulations. 
	
	\subsubsection{Simple Pore Shapes}
	
		\begin{figure}
		\begin{subfigure}{0.31\textwidth}
			\centering
			\includegraphics[width=0.90\textwidth]{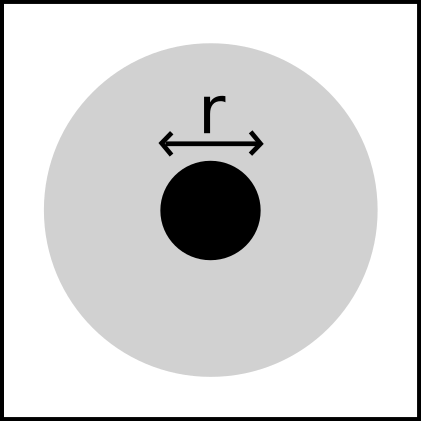}   
			\caption{}
		\end{subfigure}
		\begin{subfigure}{0.31\textwidth}
			\centering
			\includegraphics[width=0.90\textwidth]{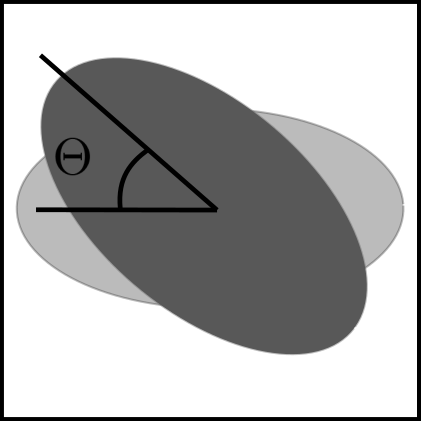}   
			\caption{}
		\end{subfigure}
		\begin{subfigure}{0.31\textwidth}
			\centering
			\includegraphics[width=0.90\textwidth]{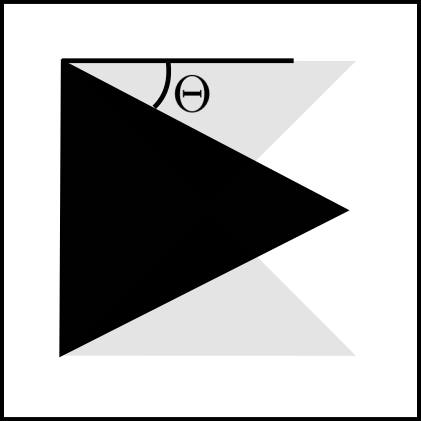}   
			\caption{}
		\end{subfigure}
		\caption{Various simple pore shape geometries: circles of varied radius, ellipses of varied pitch, and triangles of varied angle. }
		\label{fig:unitcells_simple}
	\end{figure}

	In order to evaluate simple, easily describable pore geometries, unit cells containing single shapes as a solid inclusion can be used. Although these simple unit cells are not representative of real chaotic porous media, their shapes are simple to describe and their general effect on flow and transport is easy to imagine, making the verification of results easier. Such geometries are common within the field of homogenization (\eg \cite{eggenweiler2021a}), experimental work (\eg \cite{weishaupt2020a}), or in direct numerical simulations of flow in porous media (\eg \cite{chu2019a}). 
	
	A series of simple unit cells is provided in this framework containing shapes of different types, including: squares, rectangles, circles, ellipses, triangles, and crosses. Each shape is varied in size and rotation. An example of these cells is shown in \ref{fig:unitcells_simple}. Variations and extensions of the shapes provided can be made. As these geometries do not intersect with boundaries and do not include hollow disconnected features, they do not require any preprocessing.

	\subsubsection{Randomized geometries}	
	Pseudo pore geometries are also common tools for pore-scale analysis. These geometries may not come directly from real porous media, but can be controlled and designed to reflect similar properties without requiring samples and their scans. Multiple generators exist for their creation, and a toolbox of these generators is maintained by the PoreSpy project. \cite{gostick2019a}. In this case, the voronoi polygon cells shown come directly from the generators in this toolbox, with a custom modification to the padding method to ensure periodic connectivity, as discussed in \cite{fritzen2008a}. Here, variations to the number of seeds and the aperture of the polylines can be made. Two further generators are added, also inspired by the perlin and fractal noise (\cite{perlin2002a}) generators developed in PoreSpy, but with customizations to guarantee periodic connectivity. These are implemented using a perlin noise generator tool based on NumPy \cite{vigier2018a}. Here variations to the scale and the filter parameter can be made. 
	
	\begin{figure}
		\begin{subfigure}{0.31\textwidth}
			\centering
			\includegraphics[width=0.90\textwidth]{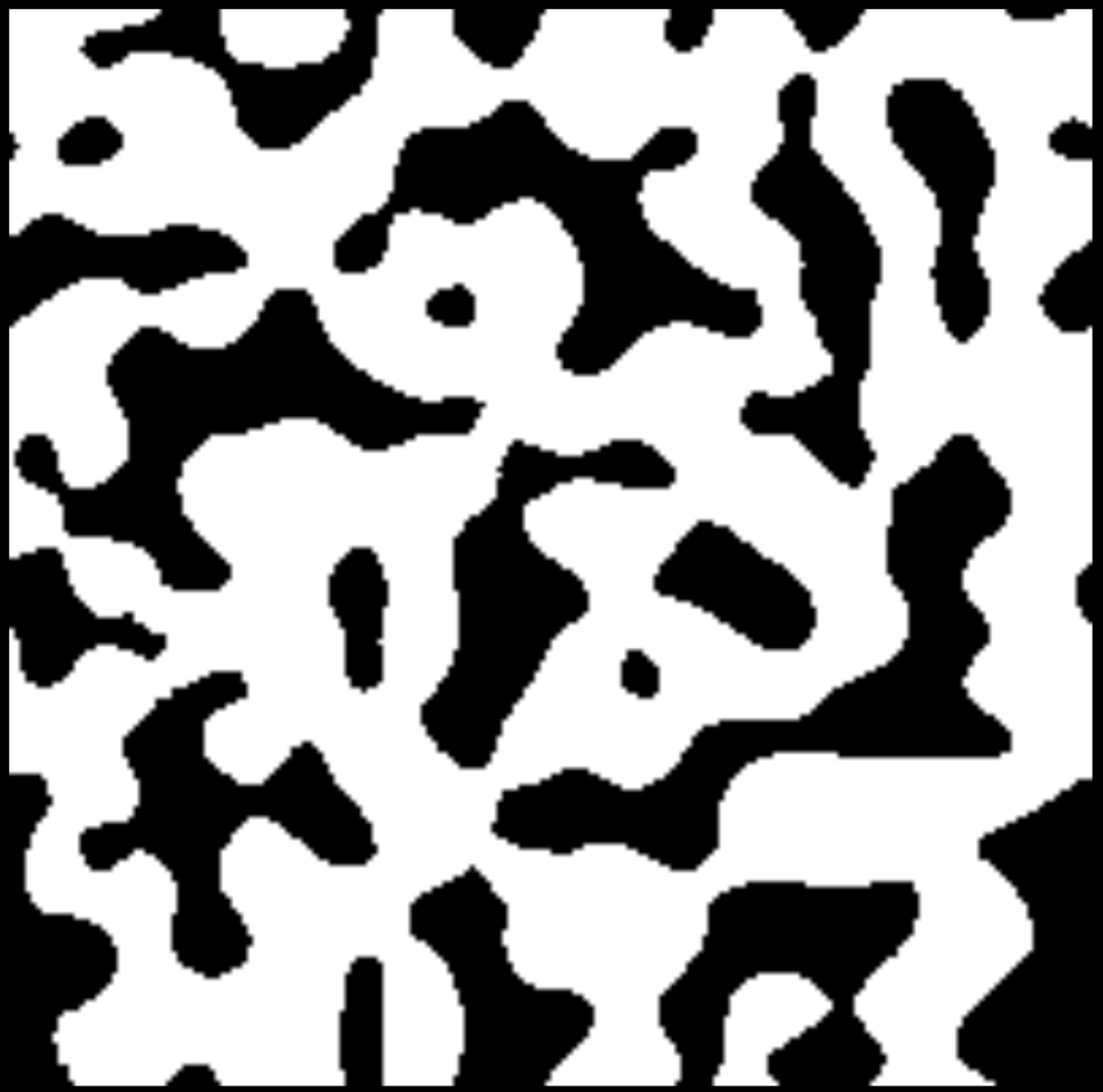}
			\caption{}
			\label{fig:unitcell_perlin}
		\end{subfigure}
		\begin{subfigure}{0.31\textwidth}
			\centering
			\includegraphics[width=0.90\textwidth]{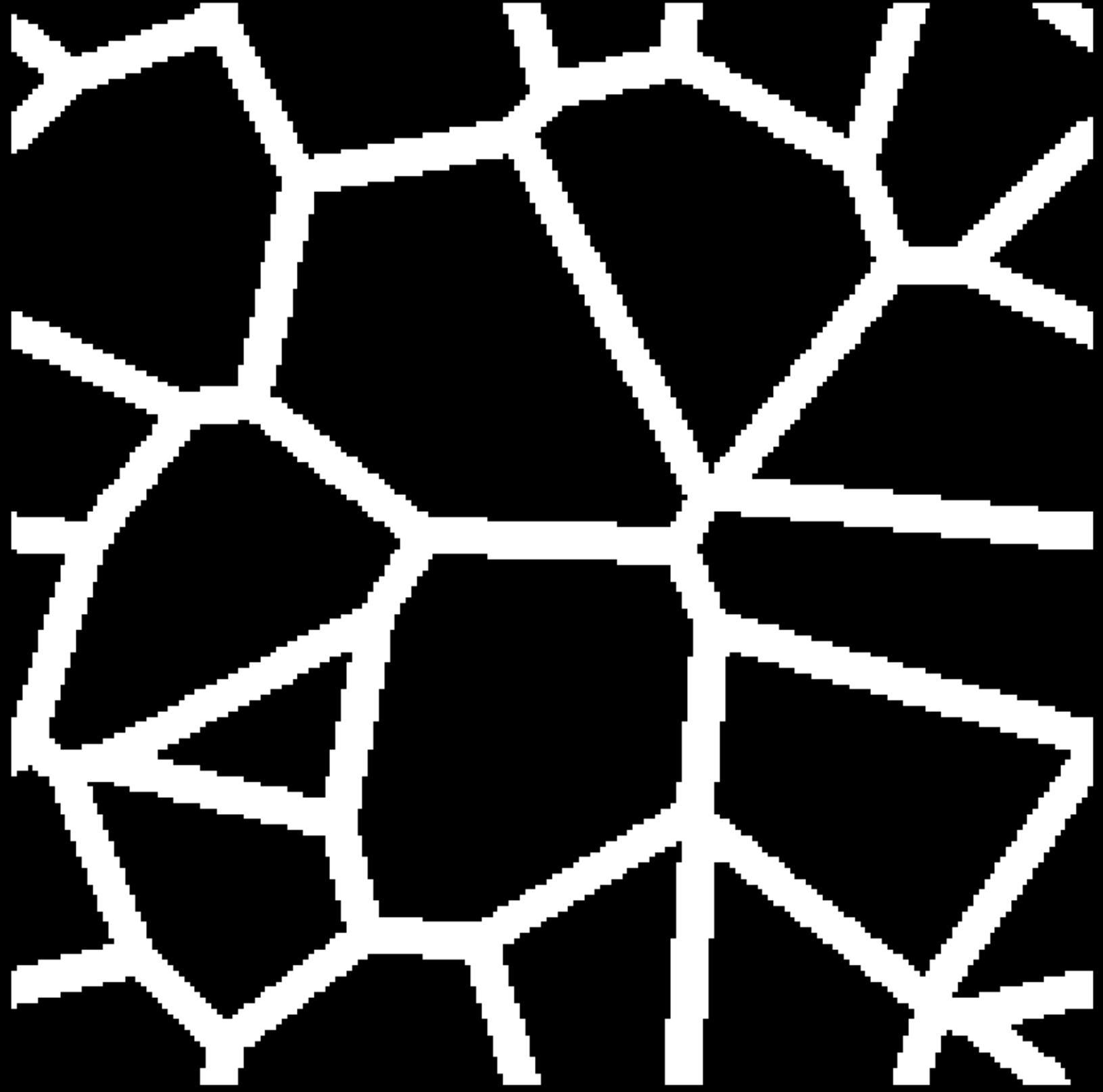}
			\caption{}
			\label{fig:unitcell_voronoi}
		\end{subfigure}
		\begin{subfigure}{0.31\textwidth}
			\centering
			\includegraphics[width=0.90\textwidth]{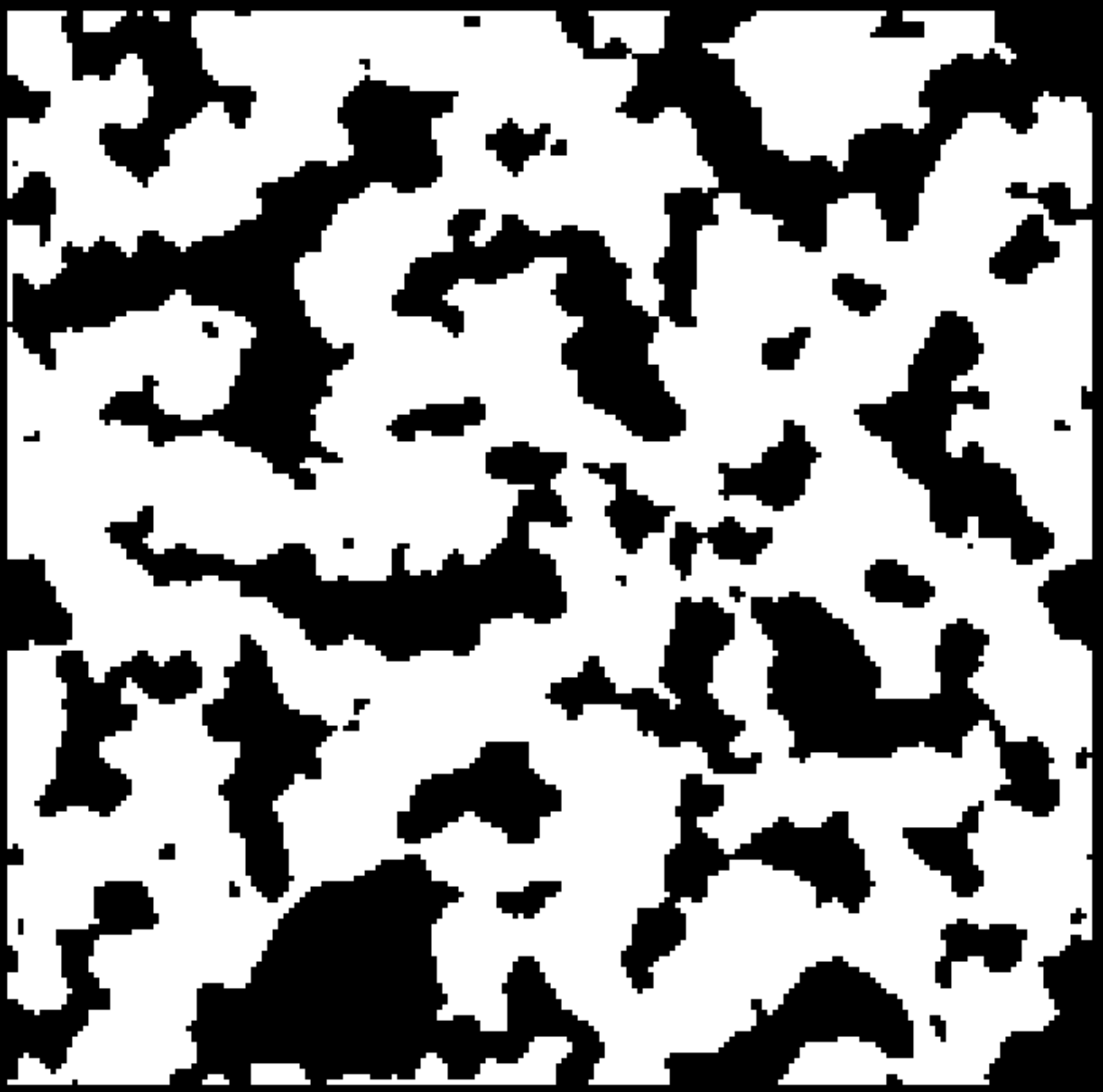}
			\caption{}
			\label{fig:unitcell_fractal}
		\end{subfigure}
		\caption{Various pseudo unit pore geometries: perlin noise based topography, voronoi polygons, and fractal noise based topography.}
		\label{fig:unitcells_complex}
	\end{figure}

	\subsection{Data Collection}
	
	When performing many simulations at the pore scale, it can become difficult to maintain an easy overview of all of the test cases and their results. In order to best curate the simulation runs in this step, \texttt{JSON} based metadata collection, implemented in \cite{lohmann2022a}, and adapted for \Dumux can be used. 
	
	In this framework a few examples of data collection in \texttt{JSON} based output are shown. Benefits of this method are the flexible type storage, and the ease of use in post-processing programs in \texttt{python} or \texttt{C++}.

%%%%%%%%%%%%%%%%%%%%%%%%%%%%%%%%%%%
%%%%%%%%%%%%%%%%%%%%%%%%%%%%%%%%%%%

	\section{Averaging and REV-Scale Model Parametrization} \label{sec:averagescale}
	
	When evaluating porous media on both the pore scale and the averaged REV scale, averages must be taken in order to compare results. Results at the pore scale will contain a high resolution solution, where results on the REV scale represent only a result based on an averaged concept \cite{quintard1993a}. Averaged results of the pore-scale model will not be equal to the the REV-scale model as some aspects of the pore-scale model will not be invariant to averaging. In order to match these variations, REV-scale models can be modified to include a closure problem, based on the underlying pore-scale problem. The development of effective parameters on the REV scale based on pore-scale problems is extensively studied in the field of Homogenization in Porous Media \cite{hornung1997a}.
	
	\subsection{Volume Averaging} \label{subsec:volume_averaging}
	
\begin{figure}
		\begin{subfigure}{0.31\textwidth}
			\centering
			\includegraphics[width=0.90\textwidth]{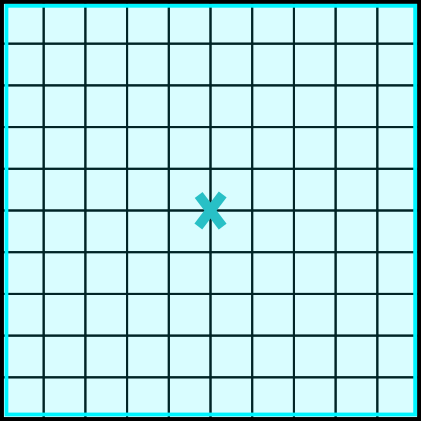}
			\caption{}
			\label{fig:volumeaveraing_full}
		\end{subfigure}
		\begin{subfigure}{0.31\textwidth}
			\centering
			\includegraphics[width=0.90\textwidth]{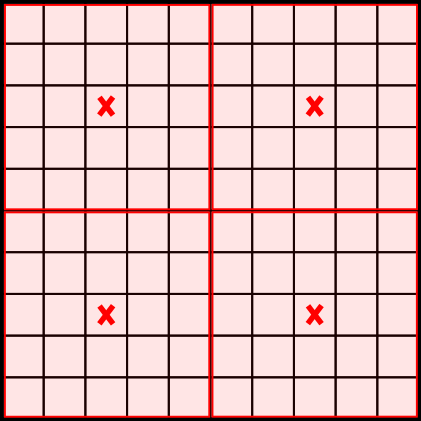}
			\caption{}
			\label{fig:volumeaveraing_sub}
		\end{subfigure}
		\begin{subfigure}{0.31\textwidth}
			\centering
			\includegraphics[width=0.90\textwidth]{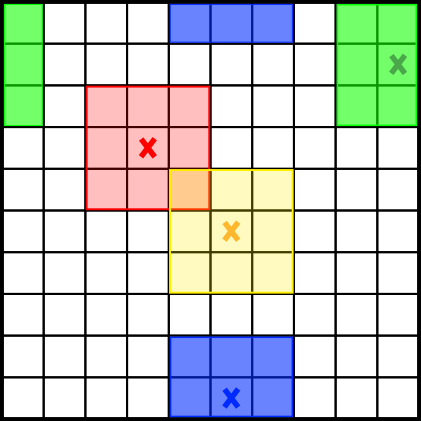}
			\caption{}
			\label{fig:volumeaveraing_conv}
		\end{subfigure}
		\caption{Various volume averaging methods: full volume averaging, sub-volume averaging, and periodic convolutional averaging.}
		\label{fig:volumeaveraing}
	\end{figure}

	In order to compare model results across scales, a spatial averaging of the subscale results to the averaged scale must be made, commonly referred  to as volume averaging \cite{whitaker1998a}. Various volume averaging methods are provided within this framework based on intersections of pore-scale and averaged scale grid geometries. In particular, three averaging methods are provided, a full volume average, a sub-averaging method, and a convolutional averaging method. Graphics depicting these methods are seen in Figure \ref{fig:volumeaveraing}. First, a full averaging method collects all pore-scale values and finds their volume-weighted average, for one unit cell, one average quantity per variable is then collected. Second, sub-averages can be taken, where the pore scale is split into multiple regions, and one average quantity is collected per region. Third, a periodic convolutional averaging can be done with a given filter or volume size. This results in a full field of averaged results at the same resolution as the pore-scale solution. 
	
	After establishing an average value for a specific scale parameter, local deviations from this average value can be defined. This can be seen with the following equation:
	
	\begin{equation}
		\underbrace{m}_{\text{Pore Scale Field}} = 
		\underbrace{\left<m\right>}_{\text{Pore Scale Field (Spatial Average)}}
		+ \underbrace{\tilde{m}}_{\text{Pore Scale Field (Variation)}}
	\end{equation}

	Pore scale field variations can typically be ignored, as their average equals zero, but in instances where these variations are in operation with other independently developed variations, a pore-scale variation term on the averaged scale will arise (\eg $\left<\tilde{m}\tilde{n}\right> \neq 0$). 

	\subsection{Closure Problem}
	
	As the averaged pore-scale model concept ($\left< \mathcal{L}_{\text{PS}}(m,n) \right>$) will contain both averaging invariant ($\left<\mathcal{L}\left(m,n\right)\right>$) and averaged variation terms ($\left<\mathcal{L}\left(\tilde{m},\tilde{n}\right)\right>$), it will not be equal to the REV-scale model concept ($\mathcal{L}_{\text{REV}}(\left<m\right>, \left<n\right> )$). The averaging invariant term ($\left<\mathcal{L}\left(m,n\right)\right>$) will be properly described by the averaged model concept ($\mathcal{L}\left(\left<m\right>,\left<n\right> \right)$), but the averaging variations are not included in the REV-scale concept. These terms, aimed to represent the effect of pore-scale variations on the REV scale, are referred to as closure problems, which can be defined using an analysis of the pore scale. These closure terms are shown here as $\mathcal{F}\left(\boldsymbol{\alpha, ...} \right)$, as can be seen here:
	
	\begin{align}
		\left< \mathcal{L}_{\text{Pore-scale}}(m,n) \right> &:=  
		\underbrace{\left<\mathcal{L}\left(m,n\right)\right>}_{\text{Averaged invariant}}
		+ \underbrace{\color{red}{} \left<\mathcal{L}\left(\tilde{m},\tilde{n}\right)\right> \color{black}{}  }_{\text{Averaged variations}} = 0, && \label{eq:transportAvg_2} \\
		\mathcal{L}_{\text{REV-scale, closed}}(\left<m\right>, \left<n\right> ) &:=  
		\underbrace{\mathcal{L}\left(\left<m\right>,\left<n\right> \right) }_{\text{Averaged model concept}} 
		+ \underbrace{\mathcal{F}\left(\boldsymbol{\alpha}, ... \right)}_{\text{Closure}} = 0. &&
		\label{eq:transportREV}
	\end{align}

	The parameters $\boldsymbol{\alpha}$ are then representative parameters based on the pore scale that can be altered on the REV scale depending on the properties of the pore scale. 
	
	\subsection{Optimization Targets}
	\label{subsec:optimization_target}
		
	With the volume-averaged results from the high-resolution pore scale, and a target closure term on the REV scale, optimal representative parameters $\boldsymbol{\alpha}$ can be found to minimze the difference between the two scales. 
	
	Using the mean least squares loss function, the following optimization target can be specified. 
		
	\begin{equation}
		\boldsymbol{\alpha} = \arg \text{min} \mathcal{L}_2 (\boldsymbol{\alpha}) , \qquad 
		\mathcal{L}_2(\boldsymbol{\alpha}) := 
		\left|\left|  \left<\mathcal{L}\left(\tilde{m},\tilde{n}\right)\right> - \mathcal{F}\left(\boldsymbol{\alpha}, ... \right) \right|\right|_2^2 
	\end{equation}

	Depending on the varied size of the parameters $\boldsymbol{\alpha}$, different loss functions can be chosen, such as the mean absolute percentace error. Optimization tools used within this framework come from the python toolkit SciPy \cite{virtanen2020a}, where various solvers are available.

	%%%%%%%%%%%%%%%%%%%%%%%%%%%%%%%%%%%
	%%%%%%%%%%%%%%%%%%%%%%%%%%%%%%%%%%%

	\section{Generalization of Parametrization with Machine Learning Tools} \label{sec:mlgeneralization}
	
	With optimal REV-scale parameters describing the solutions from pore-scale solutions, a connection between pore-scale dynamics and REV-scale parameters can be drawn. Although possible via the steps described above, it would be rather inefficient to repeatedly perform pore-scale simulations within a REV-scale model to capture better parameters. As these optimal parameters are based on the properties seen at the pore scale, a data-driven model can be developed to replace these expensive pore-scale simulation and optimization steps with a mapping between pore-scale properties and optimal REV-scale parameters. Although the pore-scale properties provided can take various forms, the form taken within this framework will be the binary image of the pore geometry. 
	
	\subsection{Convolutional Neural Networks (CNNs)}
	
\begin{figure}
		\centering
		\includegraphics[width=0.85\textwidth]{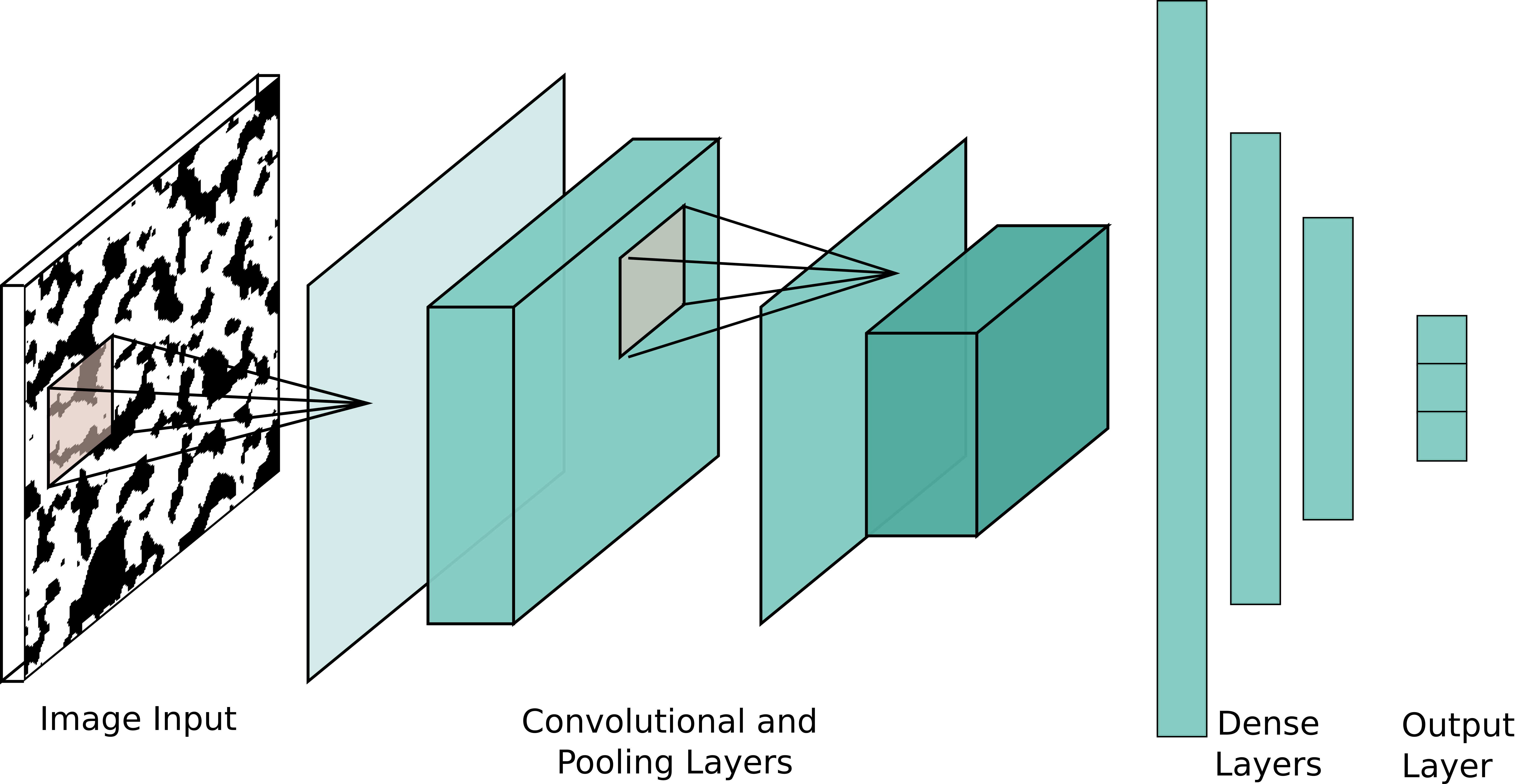}
		\caption{\label{fig:CNNs}An example visualization of a CNN (Convolutional Neural Network).}
	\end{figure}

	Artificial Neural Networks have become very powerful tools for developing data-driven mappings between input data and target solutions \cite{goodfellow2016a}. For many cases, the closure term used on the REV scale will be most effected by the form and properties of the void-solid interface in the pore scale. For this reason, a common input data form for these models will be the pore geometry \cite{gaerttner2021a, kamrava2021a}. Although standard dense layers can also be used for the processing of image input data, convolutional neural networks are effective at developing regression or classification models based on image input \cite{krizhevsky2017a}. Using a series of convolutional layers, with hidden pooling and normalization layers, followed by a series of dense layers to a regression output, this network can be trained to accurately predict the REV-scale target parameters based on the pore geometry images used as input. An example visualization of such a CNN is show in Figure \ref{fig:CNNs}
	
	All neural networks built within this framework use the toolkit TensorFlow \cite{tensorflow2015a}. Layers, connectivity, loss functions and data management tools are all available within this toolkit, and models can be built as a modular collection of these building blocks. In addition to the layers available, one additional custom layer has been built and introduced to the model to account for the periodicicty of the pore geometry images. This layer, based on the work presented in \cite{schubert2019a}, adds a periodic padding to each convolutional layer allowing wrap around filters in each direction. Training sets are initially split into three sets, a training set, a validation set, and a test set. The networks are trained on the data provided in the training set, and a stopping criteria is developed using the validation set. If, after a predefined number of training steps, the loss function as evaluated on the validation set, rather than the training set, does not decrease, the training will stop. This is a method used to avoid over-fitting. The test set can then be used to evaluate the overall accuracy of the neural network.

	\subsection{Metrics Based Neural Networks}
	
	In addition to image based input layers, other input information can be introduced to the machine learning models. For example, certain descriptive metrics can be added directly to the dense layers as additional input sources. In cases where this input layer is relevant to the target output parameters, positive effects can be seen in the training process, leading to higher learning rates and a reduced model error.

	%%%%%%%%%%%%%%%%%%%%%%%%%%%%%%%%%%%
	%%%%%%%%%%%%%%%%%%%%%%%%%%%%%%%%%%%

	\section{Analysis Using Descriptive Metrics} \label{sec:metrics}
	
	One common and debatable critique of data-driven tools in the scientific modeling context is that these models replace tools with physical meaning with black-boxes. Although this is largely accurate, in addition to their comparably computational efficiency, these models can still be used to infer about the physical processes they simulate. For example, correlation evaluations of the learned model output with descriptive metrics of the underlying system. When the learned metrics provided by the machine learning model correlate strongly with a specific metric describing one aspect of the pore-scale system, this aspect can be assumed to strongly affect the physics of the additional term. 
	
	This framework introduces a series of metrics which attempt to describe the features of the pore-scale system without solving any flow or transport specific problems. Using these metrics, correlations can be found to any predicted parameter, which can help to infer to the meaning of these REV-scale parameters. 
		
	\subsection{Porosity and Pore Size Distribution}
	\label{subsubsec:Porosity}
	
	\begin{wrapfigure}[18]{l}{0.45\textwidth}
		\centering
		\vspace{-3mm}
		\includegraphics[width=0.4\textwidth]{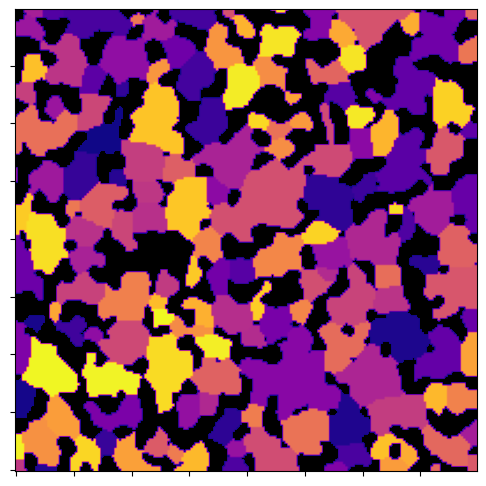}
		\caption{Example Pore Distribution Properties: $\phi = 0.65$, $\mu_{p} = 166\text{px}$, $\sigma_{p} = 140\text{px}$.}
		\label{fig:poredist}
	\end{wrapfigure}

	The most commonly used to pore geometry metrics describe the size and distribution of solid inclusions and void space within a sample. One common example of this is the porosity $\phi$, calculated as the ratio of the void space, $V_{void}$, to its total volume, $V_{total}$, $\phi = \nicefrac{V_{void}}{V_{total}}$. In addition to the porosity, a pore size distribution can classify the overall content of a pore geometry. In this case, metrics for this will include the average pore size, $\mu_{p}$, as well as the standard deviation in pore size, $\sigma_{p}$. An example pore geometry with these metrics defined is shown in Figure \ref{fig:poredist}.
	
	In this framework, the porosity is calculated as the total number of connected pore-space pixels, divided by the full number of pixels forming the image. In order to calculate the pore-size distributions, techniques from the development of pore-network models are used, specifically those outlined in detail in \cite{gostick2017a} and available in \cite{gostick2019a}. The pore space is segmented into pores by first performing a distance transform throughout all void space, applying smoothing and merging checks to mark isolated peaks, and then segmenting towards these peaks using a watershed algorithm. The result, a pore space segmented into a number of individual pores $N_{pores}$, where each has an associated pore volume. From this distribution, the statistical quantities $\mu_{p}$ and $\sigma_{p}$ can be extracted. 
	
	\subsection{Pore Surfaces: Area and Directionality}
	\label{subsubsec:SSA}
	
	As many flow conditions within a porous media are characterized by their fluid-solid interfaces, quantifications of these interfaces are useful in accurately describing pore geometries. As in some cases, these interfaces will be in contact with more than one fluid, these interfaces are described here to be gererally as void-solid interfaces. First, the specific surface area, $\textbf{S}$ ($L^{-1}$), represents the ratio between the void-solid surface area and the total volume of a porous sample \cite{bear1988a}. Second, in order to describe the general orientation of this surface area, a directionality vector can be introduced $\textbf{Di}$.

\begin{figure}
		\centering
		\begin{subfigure}{0.49\textwidth}
			\centering
			\includegraphics[width=0.80\textwidth]{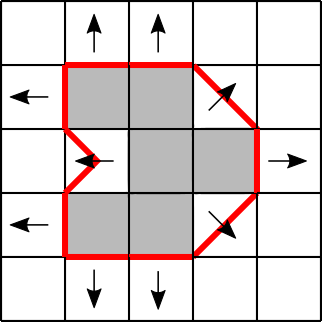}   
			\caption{An outline of the methods used for determining the surface properties. Each pixel laying on a void-solid interface contributes to the total surface area as shown in red, and the directionality of all interface adjacent pore pixels is marked with arrows.}
			\label{fig:surface_example}
		\end{subfigure}
		\begin{subfigure}{0.49\textwidth}
			\centering
			\includegraphics[width=0.80\textwidth]{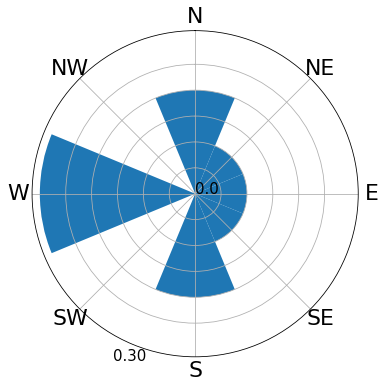}
			\caption{The directionality of the sample shown to the left visualized. The specific surface area for this sample is $\textbf{S} = 11.24$ units, and the directionality vector, and it's standard deviation are as follows: $\textbf{Di} = \left[0.2, 0.1, 0.1, 0.1, 0.2, 0.0, 0.3, 0.0\right]$, $\sigma_{\text{Di}} = 0.1$. }
			\label{fig:directionality_ex}
		\end{subfigure}
		\caption{}
		\label{fig:poro_surf}
	\end{figure}

	In this framework, the calculation of the specific surface area and the directionality vector are both calculated according to algorithms shown in Figure \ref{fig:surface_example}. First, all boundaries of pixels representing void space are evaluated. If a boundary represents a void-solid interface, this pixel contributes to the specific surface area and the directionality. If only one pixel boundary represents an interface, this intersection length is added to the surface area, and the direction of it's normal vector is added to the directionality vector. If there are two pixel boundaries adjacent to each other, the diagonal pixel length is added, and if these two boundaries are opposite each other, both interfacial lengths are added. For the directionality metric, the sum of the normal vectors is added. In the case of three pixel boundaries, a concave surface is expected, and half of a diagonal length is added twice. As in the case of two pixel boundaries, for three boundary faces, the normal vectors are collected and their sum added to the directionality vector. When evaluation of the pore-space pixels is complete, each metric is to be regularized: the specific surface area by the total area, and the directionality by the total number of boundary adjacent pixels. The results of the example case are shown in Figure \ref{fig:directionality_ex}.
	
	\subsection{Graph Based Metrics} \label{subsec:graph_methods}
			
\begin{wrapfigure}[24]{r}{0.55\textwidth}
		\centering
		\vspace{-3mm}
		\includegraphics[width=.50\textwidth]{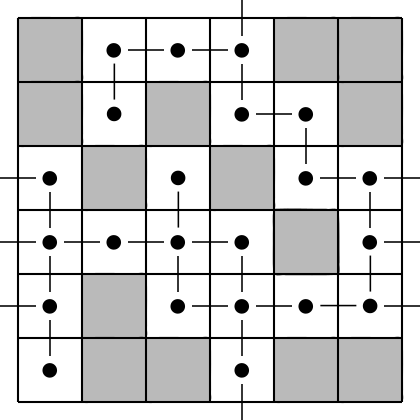}
		\caption{\label{fig:graph_scheme}An example graph representation of a pore geometry. The white and gray boxes represent the pixel-wise pore geometry. The black circles depict graph nodes, and are connected to each other via graph edges, here shown as black line segments. Edges crossing an outer boundary illustrate the periodicity of the geometry. }
	\end{wrapfigure}

	To further analyse the pore space, a metric for the flow path geometry and the general connectivity should be included to more closely describe the flow conditions in the pore space. One method for this analysis, barring the computation of any flow or transport problems, would be to represent the pore space using a graph structure. 
	Graphs, often used to describe the connectivity and structure of systems, can be built in this case from the pixel-based pore-geometry image via the following algorithm: (1) each pixel representing a void space is represented as a graph node, $N_i$, and (2) for all neighboring pore-space pixel intersections, a connecting graph edge, $E_i$, is developed connecting the adjacent nodes. Each node has a volume related to the pixel size, and each edge is non-directional and carries a uniform capacity of 1 unit. Given the periodicity of the cells evaluated in this work, boundary nodes with a direct neighbor across a periodic boundary can also be connected by an edge. The collected graph, $\text{G}$, then consists of all nodes and edges, $\text{G} = \left(\textbf{N}, \textbf{E}\right)$ \cite{harary1969a}. An example graph is shown in Figure \ref{fig:graph_scheme}.
	
	Three primary metrics can be identified using graph representations: General pore-space connectivity, a shortest path analysis to describe the geometrical tortuosity, and a maximum flow analysis. Examples of other works using graphical representations of pore geometries are as follows: \cite{kanavas2021a} where flow path resistance is investigated, \cite{tang2021a} where preferential pathways are evaluated, and \cite{gaerttner2021a} where an affinity between rock core permeability and graph's Maximum Flow is demonstrated. In addition, the development of pore-network models is a similar method, where collections of highly connected nodes form pore bodies connected by throats \cite{blunt2001a}. This graph form has also been used to relate descriptions of pore geometry to REV-scale phenomena \cite{jerauld1990a}.
		
	\subsubsection{Connectivity} \label{subsubsec:connectivity}
	
	If, from one node, a traverse can be performed such that all other nodes within the graph are reached, this graph is considered contiguous. For each additional traverse required to collect all of the pore space nodes, the connectivity increases by one unit. As is often the case in porous media, some porous regions will not be connected with the primary pore flow region. In pre-processing steps, these regions should be removed, enforcing a full connectivity. In cases where this is not performed, the degree of connectivity can be used to describe this.    
		
	\subsubsection{Shortest Path Analysis: Tortuosity} \label{subsubsec:tortuosity}

	The concept of tortuosity $\tau$ was originally introduced by \cite{carman1937a} in order to approximate the difference in flow path length between real chaotic porous media, and their common bundle-of-tubes representation. Since then, this parameter has been incorporated into many models describing flow and transport in porous materials, but has taken on multiple definitions. An in-depth discussion and summary of this topic is found in \cite{ghanbarian2013a}. Within this work, metrics used to describe pore geometries should be measureable from the pore geometry alone, meaning one should not be required to solve flow and transport problems in order to produce these simple metrics. As metrics such as the hydraulic tortuosity would require flow solutions, one is limited to a geometrical tortuosity, calculated here as: 
	
	\begin{equation}
		\tau = \frac{\left<T_{pore}\right>}{T_{L}}
	\end{equation}
	
	with $\left<T_{pore}\right>$ being the average geometrical path transversal, and $T_{L}$ being the direct straight line length. 
	
	In order to calculate the tortuous length $\left<T_{pore}\right>$, a geometric shortest path analysis is performed using the shortest path graph algorithm. From one source node, the distances from all adjacent nodes to a target node is calculated using breadth first searches. The adjacent node with the shortest distance to a target node marks the next step along the shortest path and all adjacent nodes are then evaluated again. With this, the path grows iteratively \cite{dijkstra1959a} and the distance taken along the final path can be evaluated. As the graphs produced in this work are developed from Cartesian pixels, a stair-wise traverse is favored in the implemented algorithm, and the corresponding path distance for stair-wise sub-paths is evaluated as a diagonal. To determine the tortuosity from these paths, sample locations are selected on opposite sides of the porous sample for each axis. The shortest path between each of these locations is found, and each path length is calculated. The average of these lengths is then evaluated and divided by the straight-line length. An example solution is shown in Figure \ref{fig:tortuosity}.
		
	\subsubsection{Maximum Flow} \label{subsubsec:max_flow}
	
	The Maximum Flow metric, $F_{max}$, can represent the size of the collective bottlenecks within the porous medium, and in turn, the geometric resistance to flow. Often referred to as the minimum cut, the maximum flow of a graph describes the total edge capacity connecting one source node in the graph to a sink node, or the minimum edge capacity that would disconnect the graph if removed. To evaluate this, connecting paths between the source and sink nodes are found, and the edge capacity is added to the result. The edges forming this path are then reduced by the quantity of the added edge capacity. More paths and reductions are found iteratively until the graph is no longer connected, leaving the collected capacity to be the graph's maximum flow \cite{fordfulkerson1956a}. 

\begin{figure}
		\centering
		\begin{subfigure}{0.49\textwidth}
			\centering
			\includegraphics[width=0.95\textwidth]{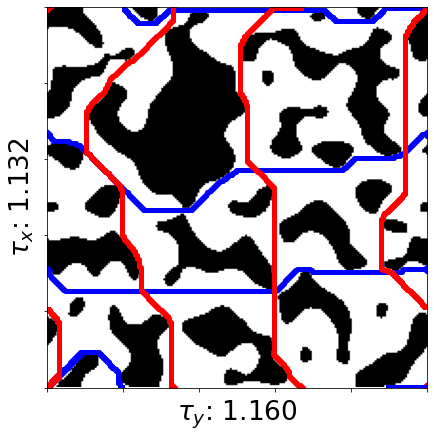}
			\caption{An example tortuosity calculation. Tortuous paths $\left<T_{pore}\right>$ are shown in blue and red for the x and y directions respectively. }
			\label{fig:tortuosity}
		\end{subfigure}
		\begin{subfigure}{0.49\textwidth}
			\centering
			\includegraphics[width=0.95\textwidth]{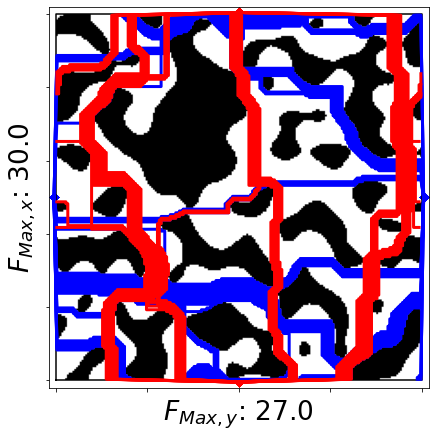}
			\caption{An example maximum flow calculation. The edges removed in order to disconnect the graphs are shown in blue and red for the x and y directions respectively.}
			\label{fig:maximumFlow}
		\end{subfigure}
		\caption{}
		\label{fig:graph_examples}
	\end{figure}

	To evaluate this for a pore geometry, a source and sink node are added to either side of a pore geometry in one direction, and periodicity is applied only to the off-axes.  Each pore geometry is analyzed in each direction to find a $F_{max}$ for each dimension. Larger differences between the axial $F_{max}$ values also indicate a degree of flow anisotropy, where similar values show a more isotropic medium. An example solution is shown in Figure \ref{fig:maximumFlow}.
	
	%%%%%%%%%%%%%%%%%%%%%%%%%%%%%%%%%%%
	%%%%%%%%%%%%%%%%%%%%%%%%%%%%%%%%%%%

	\section{Conclusions} \label{sec:conclusions}
	
	Within this work, a framework for the data-driven determination of multi-scale parametrizations is outlined. Simulations of flow and transport on the REV scale are efficient tools used ubiquitously in the environmental, biological, and technical fields surrounding porous media research. Simulations at this scale, although efficient, require parameters that represent pore-scale phenomena to perform accurately.  This work outlines a series of four steps where these parameters can be calculated from pore-scale data, their generalized and their content better understood with pore metrics. 
	
	First, a description of the pore-scale evaluation and the production of pore-scale training data is described in Section \ref{sec:porescaledata}. Second, methods for averaging and the parameterization of the REV-scale model are described in \ref{sec:averagescale}. Third, the generalization of this parameterization using a convolutional neural network is outlined in Section \ref{sec:mlgeneralization}. Finally, methods to qualify the results physically with pore scale descriptive metrics are outlined in Section \ref{sec:metrics}. 

	This framework aims to provide structure when using \Dumux to define multi-scale parameterizations. Further, the methods outlined in Section \ref{sec:porescaledata} are designed to be flexibly replaced with other pore-scale data sources, such as experimental data. Target parameterizations for this framework are focused around topics such as the dispersive transport of flow in porous media, varied momentum balances in porous-media flow models, such as Darcy-Forchheimer models, coupling conditions for coupled free flow and porous media flow models, as well as parameters required for two-phase flow models in porous media. 

\vspace{-2mm}
\section*{Acknowledgments}
\vspace{-2mm}
The authors would like to thank Prof. Dr.-Ing Andrea Beck for her guidance and constructive discussions. 
All authors would like to thank the German Research Foundation (DFG) for supporting this work by funding SimTech via Germany’s Excellence Strategy (EXC 2075 – 390740016).
The authors would also like to thank the DFG for supporting this work by funding SFB 1313, Project Number 327154368, Research Project A02.
\vspace{-5mm}
\section*{Declarations}
\vspace{-2mm}
The authors do not have any relevant financial or non-financial affiliations to declare.

\bibliography{literature}

\end{document}